\documentclass[11pt,a4paper,leqno,twoside]{amsart}
\usepackage{latexsym,amssymb,amsmath}
\input xy
\xyoption{all}
\usepackage[english]{babel}

\usepackage{amsthm}
\usepackage{amscd}
\usepackage{tabularx}

\def\o{\omega}

\def\cE{{\mathcal E}}

\def\cC{{\mathcal C}}

\def\cO{{\mathcal O}}

\def\CC{\mathbb C}

\def\ZZ{\mathbb Z}
\def\GG{\mathbb G}
\def\PP{\mathbb P}

\def\fsl{{\mathfrak {sl}}}
\def\fgl{{\mathfrak {gl}}}

\def\om{\omega}

\def\g{\gamma}

\def\ot{{\mathord{\,\otimes }\,}}
\def\op{{\mathord{\,\oplus }\,}}

\def\lra{{\mathord{\;\longrightarrow\;}}}
\def\ra{{\mathord{\;\rightarrow\;}}}
\def\we{{\mathord{{\scriptstyle \wedge}}}}

\newcommand\rem{{\medskip\noindent {\em Remark}.}\hspace{2mm}}

\newtheorem{theo}{Theorem}
\newtheorem{coro}[theo]{Corollary}
\newtheorem{lemm}[theo]{Lemma}
\newtheorem{prop}[theo]{Proposition}

\begin{document}

\title{On linear spaces of skew-symmetric matrices of constant rank}
\author{L. Manivel, E. Mezzetti}
\maketitle

\section{Introduction}

Linear sections of the Grassmannians $\GG(1,n)$ of lines in $\PP^n$
appear naturally in several different situations. In complex
projective algebraic
geometry,  3-dimensional linear sections of
$\GG(1,4)$ appear in
the classification of Fano threefolds,  2-dimensional
linear sections of
$\GG(1,5)$ define one of the smooth scrolls of $\PP^5$.

Linear sections of dimension $n-1$ of the Grassmannian of lines of $\PP^n$ are
classically called {\it linear congruences}. Recently Agafonov and
Ferapontov have
introduced a construction
     establishing an important connection between certain
  hyperbolic systems of conservation laws, called
of Temple class,  and congruences of lines
(\cite{AF}, \cite{AF1}). In particular they have proved that
linear congruences in $\PP^3$ and $\PP^4$ up to projective
equivalence correspond
bijectively  to Temple systems in 2 and 3 variables up to reciprocal
transformations.

The classification of the linear congruences in $\PP^3$ has been done
classically by
Kummer (\cite{K}) (he classified more generally congruences of order 1, his
classification has been completed by Ran \cite{R}), in
$\PP^4$ the classification was made by Castelnuovo (\cite{C}). In
higher dimension, the
classification is still missing.

\smallskip
A linear congruence $B$ in $\PP^5$ is  of the form $B=\GG(1,5)\cap
\Delta$, where $\Delta$ is a linear space of dimension 10. In the
Pl\"ucker embedding of
$\GG(1,5)$ in $\PP^{14}$ a line $\ell\in B$ is represented by a
skew-symmetric $6\times
6$ matrix
    $(p_{ij})_{i,j=0,\ldots,5}$ of rank two. $\Delta$ is the intersection
of four hyperplanes,
whose equations can be written in the form
$$ \sum_{i,j=0}^5 a_{ij}p_{ij}=0.$$
The four hyperplanes are
points of the dual space  ${\check \PP}^{14}$, they generate the 3-space
${\check\Delta}$.

The dual variety of the Grassmannian parametrizes the tangent
hyperplanes to $\GG(1,5)$. It is the cubic hypersurface
${\check\GG}(1,5)$ in ${\check \PP}^{14}$ defined by the Pfaffian
of the skew-symmetric matrix $(a_{ij})$, and parametrizing matrices 
of rank at most four.
Hence the intersection
${\check\GG}(1,5) \cap {\check \Delta}$  is, in general, a cubic surface $S$.

The rational Gauss map $\gamma: \check\GG(1,5) \dasharrow \GG(1,5)$ associates
to a tangent hyperplane its unique tangency point. $\gamma$ is
regular outside the singular
locus of $\check\GG(1,5)$, which is naturally isomorphic to $\GG(3,5)$ and
is formed by the
hyperplanes whose associated matrix  $(a_{ij})$ has rank 2. The fibre
of $\gamma$ over a
line $\ell$ is a linear space of dimension 5, naturally identified with
the linear system of
hyperplanes containing the tangent space to $\GG(1,5)$ at the point
$\ell$. Its intersection with
$\GG(3,5)$ is a smooth 4-dimensional quadric, representing the
3-spaces containg $\ell$.

For general $\Delta$, $S$ does not intersect $\GG(3,5)$ and the image
$\gamma(S)$ is a
2-dimensional family of lines, whose union is a smooth 3-fold $X$ of
degree 7, called
Palatini scroll. It results that the lines of the congruence $B$ are
the 4-secant lines of
$X$ and
$X$ is the fundamental locus of $B$ (\cite{FM}).

Classifying linear congruences in $\PP^5$ amounts to describing all
special positions of
the 3-space $\check\Delta$  with respect to $\check\GG (1,5)$ and to
its singular
locus.

\smallskip
An interesting case arises when the cubic surface $S$ (not intersecting
$\GG(3,5)$) has a plane as irreducible component.   Such a plane can
be interpreted
naturally as a  linear system of skew-symmetric matrices of constant rank 4, of
(projective) dimension 2.

The problem of the classification of linear systems of matrices
has been considered in several different contexts.
It is a classical problem in linear algebra to determine the maximal
dimension of a space of matrices (possibly symmetric or skew-symmetric)
of constant rank; upper and lower bounds
are known but the exact answer seems to be unknown in general
(see for instance \cite{IL} and references therein). A connected interesting
question would be to characterize the maximal spaces of matrices of 
constant rank.

The aim of this paper is to give a complete classification of the linear
systems of $6\times 6$ matrices of constant rank 4, or, in geometrical
terms, of linear subspaces of
$\check\PP^{14}$ contained in $\check\GG (1,5)$ and not
intersecting its singular
locus, up to the action of the projective linear group $PGL_6$.

If $\Delta$ is such a subspace,  the restriction of $\gamma$ to
$\Delta$ is regular and
defined by  homogeneous polynomials of degree 2 (the derivatives of the
Pfaffian). It must be injective (its fibers are linear spaces of zero 
dimension,
since otherwise they would intersect the singular locus of $\check\GG(1,5)$,
which cuts every fiber of $\g$ along a quadric hypersurface),
hence a double
Veronese embedding.  The image
$\gamma(\Delta)$ is therefore a Veronese variety embedded in $\GG(1,5)$.
All double Veronese
embeddings of linear spaces in Grassmannians of lines have been
recently classified
(\cite{SU}).  These embeddings are given by vector bundles $\cE$ of rank 2
with det$(\cE)=2$.

In the case of $\GG(1,5)$ there are the following possibilities:
\begin{itemize}
\item[(i)] $\PP^1$ can be embedded as  $v_2(\PP^1)$  by
$\cO_{\PP^1}(1)\oplus \cO_{\PP^1}(1)$
or by
$\cO_{\PP^1}\oplus \cO_{\PP^1}(2)$, it
represents  the lines of a ruling of a smooth 2-dimensional quadric
in the first case or the
lines of a quadric cone in the second one,

\item[(ii)] $\PP^3$ can be embedded as  $v_2(\PP^3)$ only as the
family of lines contained in a
smooth 3-dimensional quadric,  the bundle $\cE$  is the cokernel of a
general map
$\cO_{\PP^3}\rightarrow\Omega_{\PP^3}(2)$,

\item[(iii)] $\PP^2$ can be embedded as  $v_2(\PP^2)$ in 4
ways, via  the following bundles:
$\cO_{\PP^2}(1)\oplus \cO_{\PP^2}(1)$,
$\cO_{\PP^2}\oplus \cO_{\PP^2}(2)$, the restriction of the bundle
$\cE$ of the previous case, or
the Steiner bundle; the corresponding families of lines are
respectively: the lines
joining the corresponding points of 2 disjoint planes, the lines of a
cone over a Veronese
surface, the lines contained in a smooth 3-dimensional quadric and
meeting a fixed line in it, or
the secant lines of a skew cubic in a 3-space.
\end{itemize}

For each of the cases (i) and (iii) we show in \S 3 and 4 that there
is a linear space $\Delta$
with
$\gamma(\Delta)$ the corresponding Veronese variety. We show also
that each case gives a
unique orbit under the action of $PGL_6$. In fact, our main result
is that the space of projective planes of skew-symmetric matrices
of order $6$ of constant rank four, has four connected components,
each of which is a $PGL_6$-orbit of dimension $26$.
On the contrary, in the case
(iii), we show that there
are no 3-planes parametrizing matrices of constant rank four.

The methods of proof are mainly algebraic, in particular we reverse
the point of view and
interpret the spaces we are interested in as subspaces of $\PP^{14}$
rather than its dual. This
allows to interpret generators of these spaces  as tensors of tensor
rank two, and to give
explicit bases for the vector spaces we consider. Note that with this
point of view, the rational map $\g$ is particularly simple: it maps
a skew-symmetric tensor $\o\in\wedge^2\CC^6$ to $\o\we\o \in\wedge^4\CC^6\simeq
\wedge^2(\CC^6)^*$.

\smallskip
%\section*{Notation}
We shall always work over the complex field $\CC$.

\section{Planes of skew-symmetric matrices of order $5$}

We begin with the case of skew-symmetric matrices of order $5$.
Those of rank two are parametrized by the Grassmannian
$\GG(1,4)\subset\PP\wedge^2\CC^5\simeq\PP^9$. Since its
dimension equals $6$,
a $\PP^3\subset\PP^9$ will always meet $\GG(1,4)$. The set of planes
$\PP^2\subset\PP^9$
meeting $\GG(1,4)$ is an irreducible divisor in $\GG(2,\PP\wedge^2\CC^5)$,
the Chow hypersurface
of $\GG(1,4)$. The complement of this divisor is the space of
projective planes of
skew-symmetric matrices of constant rank equal to four.

If we choose a basis $e_1,e_2,e_3,e_4,e_5$ of $\CC^5$, an example is
given by the projective
plane $\pi_5$ generated by
\begin{eqnarray}\nonumber
e_1\we e_4+e_2\we e_3, \\ \nonumber
e_1\we e_5+e_2\we e_4, \\ \nonumber
e_2\we e_5+e_3\we e_4.
\end{eqnarray}

\begin{prop}
The $PGL_5$-orbit of $\pi_5$ is an open subset of
$\GG(2,\PP\wedge^2\CC^5)$, whose complement
is exactly the Chow hypersurface of $\GG(1,4)$.
\end{prop}

\proof A straightforward computation shows that the stabilizer of
$\pi_5$ in $PGL_5$
has for Lie algebra the space
$$\cC=\Big\{
\begin{pmatrix} 2x& -2z & 0 & 0 & 0\\ -y & x & z & 0 & 0 \\ 0 & 3y &
0 & 3z  & 0 \\
0 & 0 & y & -x & z \\ 0 & 0 & 0 & 2y & -2x \end{pmatrix}, x,y,z\in\CC \Big\}
\subset\fsl_5.$$
\nonumber
Since $\cC$ has dimension three, the orbit of $\pi_5$ is open in
$\GG(2,\PP\wedge^2\CC^5)$.
Since $\cC\simeq\fsl_2$ is reductive, by Matsushima's theorem
(\cite{matsu}) the
orbit of $\pi_5$
is affine, so its complement must be a divisor in
$\GG(2,\PP\wedge^2\CC^5)$. 
%Since the Picard group of the Grassmannian is free of rank one, this divisor
%must be irreducible, thus equal to the Chow hypersurface of $\GG(1,4)$.
This divisor contains the Chow hypersurface of $\GG(1,4)$, and we must 
check that it has no other irreducible component. For this we prove
that any plane $\pi$ not meeting $\GG(1,4)$ has a $PGL_5$-orbit of maximal
dimension. 

First observe that the map $\gamma: \wedge^2\CC^5\ra\wedge^4\CC^5\simeq
(\CC^5)^*$, mapping $\om$ to $\om\wedge\om$, must be injective on 
$\pi$, so that three general points will be mapped to three 
hyperplanes intersecting 
along a projective line $\ell\subset\PP^4$. Let $P$ be a 
projective plane, skew to this line. Our plane $\pi$ is 
contained in $\ell\wedge P\op\wedge^2\ell$, and its intersection 
with the hyperplane $\ell\wedge P$ must be a line (otherwise it would
meet $G(1,4)$). 
Let $e_1,e_2$ be a basis of $\ell$ and $e_3,e_4,e_5$ a basis of $P$. 
We may suppose that $\pi$ is generated by 
\begin{eqnarray}\nonumber
 & e_1\we e_2+e_3\we e_4, \\ \nonumber
 & e_1\we e_3+e_2\we f, \\ \nonumber
 & e_1\we e_5+e_2\we g,
\end{eqnarray}
where $f=f_3e_3+f_4e_4+f_5e_5$ belongs to $P$, and we may 
suppose that $g$ equals $e_3$ or $e_4$. 
(Note that the last two tensors 
generate the intersection of $\pi$ with $\ell\wedge P$.) 
A computation shows that $\pi$ has constant rank four
if and only if $g=e_3$ and $f_4,f_5\neq 0$, or 
$g=e_4$ and $f_3,f_4,f_5\neq 0$. Up to a change of basis
we may suppose in both cases that $f=e_4+f_5e_5$ with
$f_5\neq 0$. Then we compute explicitely that the Lie algebra of the 
stabilizer of $\pi$ in $PGL_5$ has dimension $3$, and we are done. 
\qed

\rem The appearance of this $\fsl_2$ can be explained as follows.
Let $U$ be a two dimensional vector space, and let us identify $\CC^5$ with
$S^4U$. Then, as a $SL(U)$-module,
$$\wedge^2(S^4U)=S^6U\op S^2U,$$
and one can check that the projectivization of this copy of $S^2U$, 
which has dimension
three, parametrizes a plane of skew-symmetric matrices of constant 
rank four. Of course
the stabilizer of this plane will contain (and in fact be equal to) a 
copy of $SL(U)$.

\section{Lines of skew-symmetric matrices of order $6$}

We turn to skew-symmetric matrices of order $6$. The orbit closures of
$PGL_6$ in $\PP\wedge^2\CC^6\simeq\PP^{14}$ are $\GG(1,5)\subset Pf
\subset\PP\wedge^2\CC^6$,
where $Pf$ denotes the cubic Pfaffian hypersurface of matrices of
rank at most four.

We are interested in projective lines of matrices of constant rank
equal to four, that is,
lines contained in $Pf$ but that do not touch $\GG(1,5)$. Examples are
given by the line $\ell_g$
generated by
\begin{eqnarray}\nonumber
e_0\we e_2+e_1\we e_3, \\ \nonumber
e_0\we e_4+e_1\we e_5,
\end{eqnarray}
    and the line $\ell_s$ generated by
\begin{eqnarray}\nonumber
    e_0\we e_2+e_1\we e_3, \\ \nonumber
    e_0\we e_4+e_1\we e_2.
\end{eqnarray}

\begin{prop}
A line $\ell$ of matrices of constant rank equal to four, is
$PGL_6$-equivalent either to
$\ell_g$ or to $\ell_s$.
\end{prop}

\proof Choose two points $\o$, $\o'$ generating $\ell$. They define
skew-symmetric endomorphisms
of $\CC^6$ whose images $L$, $L'$ have dimension $4$. They cannot be
equal, since otherwise the
line $\ell$ would be contained in $\PP\wedge^2L\simeq\PP^5$ and would
therefore necessarily meet
the quadric hypersurface $Q_L$ of rank two matrices. So the
intersection $L\cap L'$ must have
dimension $2$ or $3$.

\smallskip
Suppose the dimension is two, and let $e_0,e_1$ be a basis of $L\cap
L'$. If $\o\in\wedge^2L$
does not belong to the tangent space of $Q_L$ at $e_0\we e_1$, we can
find two other vectors
$e_2,e_3$ in $L$ such that $\o=e_0\we e_1+e_2\we e_3$ (up to scalar).
Similarly, if $\o'\in\wedge^2L'$
does not belong to the tangent space of $Q_{L'}$  at $e_0\we e_1$, we
can find two other vectors
$e_4,e_5$ in $L'$ such that %*
$\o'=e_0\we e_1+e_4\we e_5$.
But then the
line $\ell$ is not contained
in the Pfaffian hypersurface. If $\o'$ is contained in the tangent
space of $Q_{L'}$  at $e_0\we e_1$,
we can find two other vectors $e_4,e_5$ in $L'$ such that $\o'=e_0\we
e_4+e_1\we e_5$, and again
the line $\ell$ is not contained in $Pf$. We conclude that $\o$ and
$\o'$ both belong respectively
to the tangent space of $Q_L$ and $Q_{L'}$ at $e_0\we e_1$, so that
$\o = e_0\we e_2+e_1\we e_3$ %*
and $\o'= e_0\we e_4+e_1\we e_5$
for a suitable basis %*
$e_0,\ldots ,e_5$. Thus $\ell$ is
$PGL_6$-equivalent to $\ell_g$.

\smallskip
Now suppose that the dimension is three.
%, with basis $e_0,e_1,e_2$.
Choose vectors $e_3$ in $L-L'$ and $e_4$ in $L'-L$.
We can write $\o=\phi+e_1\we e_3$ and $\o'=\phi'+e_2\we e_4$, with
$e_1,e_2\in L\cap L'$
and $\phi, \phi'\in \wedge^2(L\cap L')$. Note that $e_1$ and $e_2$
must be %*
independent,
since otherwise we can find a linear combination of $\phi$ and
$\phi'$ of the form
$e_1\wedge f$, and the corresponding combination of $\o$ and $\o'$
would have rank two.
Now choose a vector $e_0\in L\cap L'$, not on the line joining $e_1$ to $e_2$.
Changing $e_3$ and $e_4$ if necessary,
we can then write, up to scalar, $\o=e_0\we e_2+e_1\we e_3$ and
$\o'=e_0\we e_1+e_2\we e_4$.
Thus $\ell$ is $PGL_6$-equivalent to $\ell_s$.
\qed

\medskip Lines equivalent to $\ell_s$ are contained in the linear
span of a copy of
$\GG(1,4)\subset \GG(1,5)$. We call them {\it special lines}.
    Lines equivalent to $\ell_g$, which we call {\it general lines},
are contained in a unique tangent space to $\GG(1,5)$. We call the
corresponding line in $\PP^5$ the {\it pivot} of the general line.

\smallskip
Computing the stabilizers of a general and of  a special line, we check that:

\begin{coro}
The space of lines of matrices of rank four is irreducible of dimension $22$,
with an open $PGL_6$-orbit of general lines, and a codimension one
orbit of special lines.
\end{coro}

\begin{rem} 1.
If we adopt the dual point of view, we have that both types of lines
have a conic as image
via the Gauss map. A conic coming from a line $\ell$ equivalent to
$\ell_s$ represents the
lines of a quadric cone; its vertex $P$ determines a subgrassmannian
in $\GG(3,5)$,
isomorphic to %*
$\GG(1,5)$, formed by the 3-spaces passing through $P$
and $\ell$ is contained
in its linear span. A conic coming from a line $\ell$ equivalent to
$\ell_g$ represents
instead  a ruling of a smooth quadric $Q$. The pivot of $\ell$ can be
interpreted, by
duality, as the 3-space generated by $Q$.

2. Let us recall that the Pfaffian hypersurface $Pf$ is the secant variety
of the Grassmannian $\GG(1,5)$. A rank two tensor $\o$ of the form 
$v_0\we v_1+v_2\we v_3$
can be interpreted as a point of the secant line of $\GG(1,5)$
joining the points $[v_0\we v_1]$ and $[v_2\we v_3]$ of $\GG(1,5)$. 
The corresponding lines
in $\PP^5$ are skew, and span precisely the 3-space $\PP(L)$. The 
pivot of $\ell$ is the
intersection of these 3-spaces as $\o$ varies in $\ell$.

\end{rem}

\section{Planes of skew-symmetric matrices of order $6$}

We now consider projective planes of matrices of constant rank equal
to four, that is,
planes contained in $Pf$ but not touching $\GG(1,5)$.

\medskip\noindent {\it Example 1}. We can choose a hyperplane
$H\simeq\CC^5$ in $\CC^6$. Then the plane
$\pi_5\subset\PP\wedge^2\CC^5\subset\PP\wedge^2\CC^6$
has constant rank four, and contains only special lines.

Conversely, the plane
$\pi_5$ determines uniquely $H$, so that its $PGL_6$-orbit $\cO_5$
has dimension $21+5=26$.

In terms of skew-symmetric matrices, $\pi_5$ is the space of matrices
$$\begin{pmatrix} 0 & 0 & 0 & a & b & 0\\
    0 & 0 & a & b  & c & 0\\ 0 & -a & 0 & c  & 0 & 0\\
-a & -b  & -c & 0 & 0 & 0\\  -b & -c & 0 & 0 & 0 & 0\\
0 & 0 & 0 & 0 & 0 & 0 \end{pmatrix}.$$

If we interpret $\pi_5$ as a plane in ${\check\GG}(1,5)$, its image
by the Gauss map $\g$ is
a projected Veronese surface. It is contained in a 4-space, i.e. the
Schubert cycle of lines
passing through a point $P$, and the lines corresponding to the
points of $\g(\pi_5)$
generate a cone of vertex $P$ over a Veronese surface. $\pi_5$ is
contained in the
inverse image, by the Gauss map, of that Schubert cycle, which is the
linear span of the
subgrassmannian of
$\GG(3,5)$ parametrizing the 3-spaces passing through $P$.

\medskip\noindent {\it Example 2}. Another example is
the plane $\pi_g$ generated by
\begin{eqnarray}\nonumber
e_0\we e_4-e_1\we e_3, \\ \nonumber
e_0\we e_5-e_2\we e_3, \\ \nonumber
e_1\we e_5-e_2\we e_4.
\end{eqnarray}

In terms of skew-symmetric matrices, $\pi_g$ is the space of matrices
$$\begin{pmatrix} 0 & 0 & 0 & 0 & a & b\\
0 & 0 & 0 & -a & 0 & c\\ 0 & 0 & 0 & -b & -c & 0\\
0 & a & b & 0 & 0 & 0\\-a & 0 & c & 0 & 0 & 0\\
-b & -c & 0 & 0 & 0 & 0 \end{pmatrix}.$$

This plane has the property that it contains only general lines. A
more intrisic way
to describe it %*
is to fix a decomposition $\CC^6=A\op B$ and an
isomorphism $u :A\ra B$.
Then $\pi_g$ is equivalent to the set of tensors of the form $x\we
u(y)-y\we u(x)$,
with $x,y\in A$. Its image by the Gauss map is the set of points of the form
$x\we y\we u(x)\we u(y)$ in $\PP\wedge^4\CC^6$. They describe the Veronese
surface of lines joining the points of $A^{\perp}\simeq B^*$, to the
corresponding
points of $B^{\perp}\simeq A^*$, the %*
correspondence being given by
the transpose map $u^t$.

A computation shows that the stabilizer of $\pi_g$  in $GL_6$
is the reductive group $(GL_3\times \CC^*)\times\ZZ_2$,
so that the orbit $\cO_g$ of $\pi_g$ %*
  has dimension $36-10=26$.

\medskip\noindent {\it Example 3}.
Our next example is a  plane $\pi_t$ contained in a tangent
space to $\GG(1,5)$, generated by
\begin{eqnarray}\nonumber
e_0\we e_2+e_1\we e_3, \\ \nonumber
e_0\we e_3+e_1\we e_4, \\ \nonumber
e_0\we e_4+e_1\we e_5.
\end{eqnarray}

In terms of skew-symmetric matrices, $\pi_t$ is the space of matrices
$$\begin{pmatrix} 0 & 0 & a & b & c & 0\\
0  & 0 & 0 & a  & b & c\\ -a & 0 & 0 & 0 & 0 & 0\\
-b & -a  & 0 & 0 & 0 & 0\\ -c & -b & 0 & 0 & 0 & 0\\
0 & -c & 0 & 0 & 0 & 0 \end{pmatrix}.$$

The special lines in that plane are parametrized by a smooth conic in
the dual plane.
More intrisically, we can choose two different hyperplanes
$A$ and $B$ in a supplement $P$ of $\langle e_0, e_1\rangle$, and an
identification
$u : A\ra B$ with no fixed line. Then $\pi_t$ is equivalent to the
plane of tensors
of the form $e_0\we x+e_1\we u(x)$, with $x\in A$. Its image via the
Gauss map is
the set of points of the form $e_0\we e_1\we x\we u(x)$.
%This can be
%seen to coincide
%with the set of secant lines to a normal rational cubic curve in $\PP^3$.

Dually, the Gauss image of $\pi_t$ is the Veronese variety, contained
in the Schubert cycle
of the lines of a 3-space $S\subset\PP^5$, parametrizing the 2-secant
lines of a skew cubic
curve $C$ contained in $S$. The inverse image via $\g$ of that
Schubert cycle is the tangent
space to $\GG(3,5)$ (naturally identified with the singular locus of
$\check\GG(1,5)$) at
the point representing $S$. The special lines contained in $\pi_t$
correspond to the cones over $C$ with vertices at the points of $C$.

The Lie algebra of the stabilizer in $GL_6$ of our plane $\pi_t$,
is the subalgebra of $\fgl_6$ of matrices of the form
$$\begin{pmatrix} x+u & z & p & q & r & s \\
y & v+u & q & r & s & t \\ 0 & 0 & 3v & -y & 0 & 0 \\
0 & 0 & -3z & x+2v & -2y & 0 \\
0 & 0 & 0 & -2z & 2x+v & -3y \\
0 & 0 & 0 & 0 & -z & 3x
\end{pmatrix}.$$
In particular the dimension of the stabilizer is $10$, and the orbit $\cO_t$
of $\pi_t$ has dimension $36-10=26$. One can be slightly more
precise: the stabilizer
of %*
$\pi_t$ in $GL_6$ is a semi-direct product
$(GL_2\times\GG_m)\rtimes U$, where
$U\simeq S^4V^*\ot\det V$ as a $GL_2=GL(V)$-module, and the
multiplicative group
$\GG_m$ acts on $U$ through its tautological character.

\medskip\noindent {\it Example 4}.
Finally, we have an example of a plane $\pi_p$ containing only a
pencil of special lines, generated by
\begin{eqnarray}\nonumber
e_0\we e_3+e_1\we e_2, \\ \nonumber
e_0\we e_4+e_2\we e_3, \\ \nonumber
e_0\we e_5+e_1\we e_3.
\end{eqnarray}

In terms of skew-symmetric matrices, $\pi_p$ is the space of matrices
$$\begin{pmatrix} 0 & 0 & 0 &a & b & c\\
0  & 0 & a & c  & 0 & 0\\ 0 & -a & 0 & b & 0 & 0\\
-a & -c  & -b & 0 & 0 & 0\\   -b & 0 & 0 & 0 & 0 & 0\\
-c & 0 & 0 & 0  & 0 \end{pmatrix}.$$

The image of this plane, seen in the dual space, by the Gauss map is
a Veronese surface of
lines  contained in a smooth 3-dimensional quadric and meeting a
fixed line $r$ in that
quadric. Also this Veronese surface is contained in a 4-dimensional
Schubert cycle,
precisely in that of the lines of a hyperplane %*
meeting a fixed line. The special lines in $\pi_p$ correspond to the lines in
the quadric passing
through a fixed point of $r$.

We compute that the Lie algebra of the stabilizer in $GL_6$ of that
plane $\pi_p$
is the subalgebra of $\fgl_6$ of matrices of the form
$$\begin{pmatrix} u_{00} &u_{10} &u_{20} &u_{30} &u_{40} &u_{50} \\
0  & u_{11} &u_{21} &u_{20} & 0 & u_{30}  \\
0  & u_{12} &u_{22} &-u_{10} & u_{30} & 0 \\
0  & 0 & 0 & u_{33} & -u_{20} & -2u_{10}  \\
0  & 0 & 0 &0 &u_{44} & u_{12} \\
0  & 0& 0 & 0 &u_{21} & u_{55}
\end{pmatrix},$$
with the %*
dependence relations
\begin{eqnarray}\nonumber
u_{33} & = -u_{00}+u_{11}+u_{22}, \\ \nonumber
u_{44} & = -2u_{00}+u_{11}+2u_{22}, \\ \nonumber
u_{55} & = -2u_{00}+2u_{11}+u_{22}.
\end{eqnarray}
In particular the dimension of the stabilizer is $10$, and the orbit $\cO_p$
of $\pi_p$ has dimension $36-10=26$.

\medskip
Clearly, these four examples of planes of constant rank four belong
to distinct $PGL_6$-orbits.

\section{The main result}

\begin{theo}
A plane $\pi$ of matrices of constant rank equal to four, is
$PGL_6$-equivalent either to
$\pi_g$, $\pi_t$, $\pi_p$ or $\pi_5$.
\end{theo}

\begin{coro} The space of planes of matrices of constant rank four has
four connected components, each of which is a $PGL_6$-orbit of dimension $26$.
\end{coro}

\proof The plan of the proof is the following. If the plane $\pi$ is
contained in
the span of a $\GG(1,4)$, then we already know by Proposition 2 that it
is equivalent
to $\pi_5$. Otherwise it must contain
general lines. If we can find three general lines in general
position, then $\pi$ is
equivalent to $\pi_g$. Otherwise, all the pivots of the general lines
must pass through
a fixed point. If the pivots are not constant, $\pi$ is equivalent to
$\pi_p$. If the
pivot is constant, $\pi$ is equivalent to $\pi_t$.

\medskip\noindent {\sl First step}. We prove our first claim:

\begin{lemm} A plane $\pi$ containing only special lines must be
contained in the
span of a $\GG(1,4)$. \end{lemm}

\proof Each special line generates a hyperplane of $\CC^6$ and we
must prove that
this hyperplane is constant. Suppose this is not the case, and choose
two lines $\ell$
and $\ell'$ in $\pi$, generating dinstinct hyperplanes $H$ and $H'$.
Let $L=H\cap H'$.
We consider the line in $\pi$ joining a point $\o\in\ell-\ell'$ to a
point $\o'\in\ell'-\ell$.
Their images must meet along a three dimension space
$M_{\o,\o'}\subset L$. But the image of
a generic $\o\in\ell$ meets $L$ only in dimension two, a contradiction. \qed

\medskip\noindent {\sl Second step}. Now we can suppose that $\pi$
contains a general line,
so that in fact the generic line in $\pi$ is general. We choose three
general lines
$\ell, \ell'$, $\ell''$, and denote their pivots by $d, d', d''$.

\begin{lemm} \label{lemma} %*
Let $\o''=\ell\cap\ell'$. Suppose that $d$ and $d'$ do not meet in
$\PP^5$. Then we can
find basis $e_0, e_1$ of $d$ and $e_2, e_3$ of $d'$ such that
$\o''=e_0\we e_2+e_1\we e_3$.
\end{lemm}

\proof Since %*
$\o''$ belongs to the line $\ell$ whose pivot is $d$,  we
have that $\o\in d\wedge V$
(where $V=\CC^6$). Similary $\o\in d'\wedge V$, and since $d$ and
$d'$ do not meet, $(d\wedge V)
\cap (d'\wedge V)=d\wedge d'$. Hence the claim. \qed

\medskip\noindent Suppose that the pivots $d, d', d''$ are lines in
general position in $\PP^5$.
By %*
Lemma \ref{lemma}, we can find basis $e_0, e_1$ of $d$, $e_2, e_3$ of $d'$
and $e_4, e_5$ of $d''$
such that the intersection points $\o''=\ell\cap\ell'$,
$\o'=\ell\cap\ell''$ and $\o=\ell'\cap%*
\ell''$
can be written
\begin{eqnarray}\nonumber
    &\o'' = e_0\we e_2+e_1\we e_3, \\ \nonumber
    &\o'= e_0\we e_4+e_1\we e_5, \\ \nonumber
    &\o = e_2\we (xe_4+ye_5)+e_3\we (te_4+ze_5).
\end{eqnarray}
A computation shows that this plane is contained in the Pfaffian hypersurface,
if and only if $t=y$. But we still have some freedom on the choice of
our basis.
If we change the basis $e_0,e_1$ by some matrix $P^{-1}\in GL_2$,
then $\o''$ and
$\o'$ will be %*
preserved if we also  change  $e_2,e_3$ and $e_4,e_5$ by
the transpose
matrix $P^t$. We then get $\o = e_2\we (x'e_4+y'e_5)+e_3\we (y'e_4+z'e_5)$ with
$$\begin{pmatrix} x' & y' \\ y' & z' \end{pmatrix}=
P^t\begin{pmatrix} x & y \\ y & z \end{pmatrix}P.$$
We can therefore choose $P$ in order to obtain the matrix of the
hyperbolic quadratic
form,  and we easily conclude that $\pi$ is equivalent to $\pi_g$.

\medskip\noindent {\sl Third step}. Now we suppose that we cannot
find three lines
in $\pi$ whose pivots are in general position in $\PP^5$. The
following lemma is
certainly well known.

\begin{lemm}\label{lemma2}
Consider a family of lines in $\PP^5$, such that no three of them are
in general
position. Then either they are contained in a fixed hyperplane, or
they all meet
a fixed line.
\end{lemm}

\proof If two general lines in the family meet, they are all
contained in the same
plane or they all pass through the same point. Otherwise, choose a
generic line $\ell$ in
the family, and consider the projection from that line to some
$\PP^3$. Then any two
projected lines meet, so again %*
either they are all contained in a fixed plane, and the
original lines are all contained in a fixed hyperplane, or they all
pass through
some fixed point, and the original lines all touch the same plane
$P=P_{\ell}\subset\PP^5$,
containing $\ell$.

Now consider another general line $\ell'$ and the associated plane
$P'$. Note that
$\ell\cap\ell'=\emptyset$, but since $P\cap\ell'$ and  $P'\cap\ell$
are non empty subsets
of $P\cap P'$, these two planes must meet along a line.

Finally, any general line
will meet both $P$ and $P'$. If that's at two distinct points, this
general line will
be contained inside $P+P'$, hence in a hyperplane. Otherwise, it
meets the fixed line
$P\cap P'$. \qed

\medskip We apply %*
Lemma \ref{lemma2} to the pivots of the lines in $\pi$. We
have two cases, and in
both cases we will conclude that the pivots must pass through a fixed point.

Indeed, two general pivots cannot be
in general position, for otherwise we could choose three generic lines %*
$\ell, \ell',\ell''$ in $\pi$
whose pivots do not meet pairwise. Then by Lemma %*
\ref{lemma}, the intersection
points of %*
the
lines $\ell, \ell',\ell''$ belong to $\PP\wedge^2H$, where $H$ 
denotes the linear span of
the three pivots.
Then $\PP\wedge^2H$ would contain $\pi$, a contradiction.

So the pivots of any two lines in $\pi$ must meet. This implies
either that they are all
contained in the same plane, a possibility that we have already
excluded, or that
they all pass through a same point.

\medskip\noindent {\sl Fourth step}. The next two lemmas will
conclude the proof of the Theorem.

\begin{lemm}
If the pivots of the lines in $\pi$ all pass though a given point but
are not constant, then
$\pi$ is equivalent to $\pi_p$.
\end{lemm}

\proof Suppose that the pivot is not constant. Choose a general line
with pivot $\overline{e_0e_1}$,
joining the points $\o=e_0\we e_2+e_1\we e_3$ to $\o'=e_0\we
e_4+e_1\we e_5$, for a suitable
basis of $\CC^6$.  Choose generic lines through $\o$ and $\o'$,
respectively, with pivots
$\overline{e_0f}$ and $\overline{e_0f'}$. Let $\o''$ be their
intersection point.

We may suppose that $f$ and $f'$ belong to the span of $e_1,\ldots
,e_5$. Since $\o$ belongs to
$e_0\wedge V+f\wedge V$, we conclude that $e_1\we e_3$ belongs to
$f\we V$, so that $f$ is a
combination of $e_1$ and $e_3$. Since it is not equal to $e_1$, we
may suppose after a change of basis,
that it is equal to $e_3$. Similarly we may suppose that $f'=e_5$.
Then $\o''$ belongs to the
intersection of $e_0\wedge V+e_3\wedge V$ and $e_0\wedge V+e_5\wedge
V$, which is
$e_0\wedge V+\CC e_3\wedge e_5$. So we may write $\o''=e_0\we
h+e_3\we e_5$ for some
vector $h$ in the span of $e_1,\ldots ,e_5$. After a change of basis,
we can thus suppose that
$\pi$ is generated by $e_0\we e_2+\phi, e_0\we e_4+\phi', e_0\we
e_1+\phi''$, where $\phi,
\phi', \phi''$ belong to $\wedge^2\langle e_1,e_3,e_5\rangle$. Any
skew-symmetric form in three
dimensions is decomposable, so we may suppose that $\phi''=e_3\we
e_5$. Then, after substracting
a suitable multiple of $\o''$ to $\o$ and $\o'$, we may suppose that
$\phi =e_1\we k$
and $\phi'=e_1\we k'$ for some vectors $k,k'$ in the span of
$e_3,e_5$. These vectors
must be %*
independent, otherwise a combination of $\o$ and $\o'$ would
drop rank. In particular
$e_3\we e_5$ and $k\we k'$ are proportional, so after a new change of
basis we may suppose
that $\pi$ is generated by $e_0\we e_2+e_1\we e_3, e_0\we e_4+e_1\we
e_5, e_0\we e_1+e_3\we e_5$.
In other words, $\pi$ is equivalent to $\pi_p$.\qed

\begin{lemm}
If the lines in $\pi$  %*
all  have the same pivot, then
$\pi$ is equivalent to $\pi_t$.
\end{lemm}

\proof $\pi$ is generated by $\o=e_0\we f_2+e_1\we g_2$,  $\o'=e_0\we
f_3+e_1\we g_3$, %*
$\o''=e_0\we f_4+e_1\we g_4$, where $f_2,f_3,f_4$ and $g_2,g_3,g_4$
must be independent.
Modulo $e_0$ and $e_1$, they generate two spaces of dimension three
which must meet in
dimension two. We choose a basis $e_3$, $e_4$ of the intersection and
complete into basis
of the two spaces by $e_2$ and $e_5$, respectively. Then %*
$\pi$ is generated by
$xe_0\we e_1+e_0\we e_2+e_1\we e_3$,  $ye_0\we e_1+e_0\we e_3+e_1\we e_4$ and
$ze_0\we e_1+e_0\we e_4+e_1\we e_5$ for some scalars $x,y,z$. But we %*
can let $x,y,z=0$ by
adding a multiple of $e_0$ to $e_3$, $e_4$ and $e_5$. We conclude
that $\pi$ is equivalent to
$\pi_t$. \qed

\begin{coro}
There exists no $\PP^3$ of skew-symmetric matrices of order six and
constant rank four.
\end{coro}

\proof If there exists such a $\PP^3$, its image via the Gauss map
has to be a 2-uple embedding
of $\PP^3$ in $\GG(1,5)$. The only such embedding represents the line
contined in a smooth quadric
in $\PP^4$. So the planes in our $\PP^3$ must be all of type $\pi_p$,
and we may suppose that
it contains the plane generated by $\o=e_0\we e_1+e_3\we e_5$,
$\o'=e_0\we e_2+e_3\we e_4$,
$\o''=e_0\we e_3+e_4\we e_5$. Let $\Omega$ be another point in our
$\PP^3$. Substracting if
necessary a combination of $\o, \o', \o''$, we may suppose that
$$\Omega =\Omega_{01}e_0\we e_1+\Omega_{02}e_0\we e_2+\Omega_{12}e_1\we e_2
+\sum_{i\le 2<j}\Omega_{ij}e_i\we e_j.$$
Now we need that $\Omega +x\o+x'\o'+x''\o''$ has rank at most four
for all $x,x',x''$,
hence in particular
$\Omega\we \o\we \o =\Omega\we \o'\we \o' =\Omega\we \o''\we \o'' = 0$, and
$\Omega\we \o\we \o' =\Omega\we \o\we \o'' =\Omega\we \o'\we \o'' = 0$, and
$\Omega\we \Omega\we \o'' = 0$.
The first three conditions give $\Omega_{15}=\Omega_{24}=\Omega_{12}=0$, the
following three $\Omega_{14}+\Omega_{25}=\Omega_{13}=\Omega_{23}=0$. The last
condition is
$\Omega_{14}\Omega_{25}+\Omega_{01}\Omega_{23}-\Omega_{02}\Omega_{13}=0$,
and we deduce that  $\Omega_{14}=\Omega_{25}=0$. Finally,
$\Omega =\Omega_{01}e_0\we e_1+\Omega_{02}e_0\we e_2$ has rank two, a
contradiction. \qed

\section{The associated vector bundles}

A space $\PP^m$ of matrices of order $n$ of constant rank $r$ defines
a morphism
$$V^*\ot \cO_{\PP^m}\lra V\ot \cO_{\PP^m},$$
where $V=\CC^n$. The kernel $K$ and the image  $E$ of this
morphism are vector bundles of rank $n-r$ and $r$, respectively.
For a space of skew-symmetric matrices, $r=2s$ and $E\simeq E^*(1)$,
which implies that the splitting type of $E$ is $E_{|\ell}=\cO_{\ell}^{\op s}
\op \cO_{\ell}(1)^{\op s}$, where $\ell$ is any line. In particular, $E$ is
uniform.

\smallskip We identify the vector bundles $E$ and $K$, whose rank are
$4$ and $2$ respectively,
for each of our families of planes of skew-symmetric matrices of order
$6$ and constant rank $4$.

\medskip\noindent {\it Type $\pi_g$}. We can find a decomposition
$V=A\op B$ and an isomorphism
$u :A\ra B$ such that $\pi_g$ is the space of 2-forms of type $x\we
u(y)-y\we u(x)$, $x,y \in A$.
The image of this morphism is the direct sum of the planes $\langle
x,y\rangle\subset A$ and
$\langle u(x),u(y)\rangle\subset B$, hence
$$E=Q\op Q, \qquad K=\cO_{\PP^2}(-1)\op \cO_{\PP^2}(-1),$$
where $Q$ denotes the rank two tautological quotient bundle on $\PP^2$.

\medskip\noindent {\it Type $\pi_t$}. Our plane is the space of
morphisms of the form
$e_0\we x+e_1\we u(x)$, where we decompose a supplement to $\langle
e_0,e_1\rangle$ into the sum
of two hyperplanes $A$ and $B$, and $u:A\ra B$ is an isomorphism. The
image of this morphism
is  the sum of the planes $\langle e_0,e_1\rangle$ and $\langle
x,u(x)\rangle$, hence
$$E=\cO_{\PP^2}^{\op 2}\op \cO_{\PP^2}(1)^{\op 2},$$
while $K$ is neither split nor a twist of the quotient bundle, since
$c_1(K)=-2$
and $c_2(K)=3$.

\medskip\noindent {\it Type $\pi_p$}. Here the image of a morphism in
our plane is the
direct sum of $\langle e_0\rangle$, of a line in $\langle e_1,e_2,e_3\rangle$
and  a plane in $\langle e_1,e_2,e_3\rangle$. We deduce that
$$E=\cO_{\PP^2}\op \cO_{\PP^2}(1)\op Q,$$
and again $K$ is neither split nor a twist of the quotient bundle,
since $c_1(K)=-2$
and $c_2(K)=2$.

\medskip\noindent {\it Type $\pi_5$}. Since we are in a five dimensional space,
we get that the kernel $K=\cO_{\PP^2}\op \cO_{\PP^2}(-2)$ splits,
while $E$ fits into an
exact sequence
$$0\lra\cO_{\PP^2}(-2)\lra \cO_{\PP^2}^{\op 5}\lra E\lra 0.$$
Since $c_1(E)=2$ and $c_2(E)=4$, the bundle $E$ cannot decompose into
a direct sum
of twists of the trivial and quotient bundles.
\medskip

\providecommand{\bysame}{\leavevmode\hbox to3em{\hrulefill}\thinspace}
\bigskip

{\sc Laurent Manivel},

Institut Fourier, UMR 5582 (UJF-CNRS),

BP 74, 38402 St Martin d'H\`eres Cedex, France.

E-mail : {\tt Laurent.Manivel@ujf-grenoble.fr}

\medskip
{\sc Emilia Mezzetti},

Dipartimento di Matematica e Informatica,

Universit\`a di Trieste,

Via Valerio 12/1, 34127 Trieste, Italy.

E-mail : {\tt mezzette@univ.trieste.it}

\end{document}